\documentclass[11pt]{article}
\usepackage{amsfonts}
\usepackage{color}
\usepackage{amsmath,amssymb,comment}
\newtheorem{theo}{Theorem}[section]
\newtheorem{lem}[theo]{Lemma}

\newcommand{\mysection}[1]{\section{#1} \setcounter{equation}{0}}
\newcommand{\proof}{{\sc Proof.} \quad}
\newcommand{\proofc}{{\sc Proof} \ }
\newcommand{\be}{\begin{equation} \label}
\newcommand{\ee}{\end{equation}}
\newcommand{\bea}{\begin{eqnarray}\label}
\newcommand{\eea}{\end{eqnarray}}
\newcommand{\bas}{\begin{eqnarray*}}
\newcommand{\eas}{\end{eqnarray*}}
\newcommand{\bit}{\begin{itemize}}
\newcommand{\eit}{\end{itemize}}
\newcommand{\qed}{\hfill$\Box$ \vskip.2cm}
\newcommand{\nn}{\nonumber}
\newcommand{\R}{\mathbb{R}}
\newcommand{\N}{\mathbb{N}}
\newcommand{\pO}{\partial\Omega}

\newcommand{\hra}{\hookrightarrow}
\newcommand{\io}{\int_\Omega}
\newcommand{\na}{\nabla}
\newcommand{\Del}{\Delta}

\newcommand{\al}{\alpha}

\newcommand{\sig}{\sigma}
\newcommand{\pa}{\partial}
\newcommand{\bom}{\overline{\Omega}}
\newcommand{\Om}{\Omega}

\newcommand{\ov}{\overline}
\newcommand{\un}{\underline}
\newcommand{\wh}{\widehat}
\newcommand{\wt}{\widetilde}

\newcommand{\hs}{\hspace*}
\newcommand{\sm}{\setminus}

\newcommand{\vp}{\varphi}
\newcommand{\lbal}{\left\{ \begin{array}{l}}
\newcommand{\lball}{\left\{ \begin{array}{ll}}
\newcommand{\ear}{\end{array} \right.}

\newcommand{\abs}{\\[5pt]}

\newcommand{\adb}{\allowdisplaybreaks}
\newcommand{\tm}{T_{max}}

\newcommand{\kD}{k_D}
\newcommand{\KD}{K_D}
\newcommand{\on}{\overline{u}}
\newcommand{\hnz}{\wh{u}_0}

\begin{document}
\adb
\title{
A dimension-independent critical exponent in a nutrient taxis system
}
\author{
Michael Winkler\footnote{michael.winkler@math.uni-paderborn.de}\\
{\small Universit\"at Paderborn, Institut f\"ur Mathematik}\\
{\small 33098 Paderborn, Germany} }
\date{}
\maketitle
\begin{abstract}
\noindent 
In a ball $\Om\subset\R^n$ with arbitrary $n\ge 1$,
the chemotaxis-consumption system
\bas
	\lbal
	u_t = \na \cdot \big(D(u)\na u\big) - \na \cdot (u\na v), \\[1mm]
	0 = \Delta v - uv,
	\ear
\eas
is considered under no-flux boundary conditions for $u$, and for prescribed constant positive boundary data for $v$.
Under the assumption that $D\in C^3([0,\infty))$ satisfies
\bas
	D(\xi)\ge \kD (\xi+1)^{-\al}
	\qquad \mbox{for all } \xi\ge 0
	\qquad \qquad
	(\star)
\eas
with some $\al<1$ and some $\kD>0$, it is shown that for each nonnegative and radially symmetric 
$u_0\in \bigcup_{q>\max\{2,n\}} W^{1,q}(\Om)$, a uniquely determined global bounded classical solution exists.\abs
This complements a previous result according to which given any positive $D\in C^3([0,\infty))$ fulfilling
$D(\xi) \le \KD (\xi+1)^{-\al}$ with some $\al>1$ and $\KD>0$, one can find nonnegative radial initial data
$u_0\in C_0^\infty(\Om)$ such that no global solution exists.\abs
\noindent {\bf Key words:} chemotaxis; global existence; boundedness; critical exponent\\
{\bf MSC 2020:} 35B45, 35K51, 35Q92, 92C17
\end{abstract}
\newpage
\section{Introduction}\label{intro}
Cross-diffusion systems of the form
\be{00}
	\lbal
	u_t=\na\cdot (D(u)\na u) - \na\cdot (u\na v), \\[1mm]
	\tau v_t = \Del v - uv,
	\ear
\ee
are at the core of various models for the interaction of taxis-type migration with attractant consumption in population dynamics
(\cite{goldstein}, \cite{woodward1995}, \cite{tyson1999}, \cite{junping_shi_JDE2016}).
Despite accordingly intense activity in the literature, some fundamental questions to date seem to have remained unanswered:
In the most prototypical case when $D\equiv 1$ and $\tau=1$ and the system is considered fully isolated in the sense
that no-flux boundary conditions for both components are imposed, for instance, associated initial-boundary value problems
in $n$-dimensional domains have been found to admit globally bounded classical solutions for initial data of arbitrary size 
when $n=2$ (\cite{taowin_consumption}, \cite{jie_jiang}); for three- and higher-dimensional relatives,
only certain global weak solutions have been constructed (\cite{taowin_consumption}, \cite{lankeit_win2025}), 
while global smooth solutions are known to exist only under appropriate smallness conditions on the initial data
(\cite{tao_consumption}, \cite{heihoff_PAMS}, \cite{baghaei}, \cite{yang_ahn}).
For relatives involving different types of boundary conditions, some results going into similar directions can be found in 
\cite{FLM} and \cite{lankeit_win_NON}.
Quasilinear versions of (\ref{00}) in which still $\tau=1$, and in which $D$ suitably generalizes the prototype given by
\be{proto}
	D(\xi)=\KD (\xi+1)^{-\alpha},
	\qquad \xi\ge 0,
\ee
have been studied in \cite{cao_ishida}, \cite{chunhua_jin2} and \cite{win_ANS}, resulting in statements on global existence of
bounded solutions under the assumption that $n\in\{2,3\}$ and 
\be{01}
	\al<0.
\ee
In the presence of Dirichlet or Robin type boundary conditions for the second component, corresponding results so far
require slightly stronger restrictions on $\al$ (\cite{chunhua_jin}, \cite{tian_zhaoyin1}, \cite{tian_zhaoyin2}).\abs
With regard to the question how far this condition (\ref{01}) is optimal for such types of conclusions, only little seems known.
In fact, it appears to have been undecided for a considerably long time whether problems of the form (\ref{00}) at all are
able to generate spontaneous singularity formation, 
or whether the absorptive character of the signal evolution mechanism expressed therein rather prevents any such phenomenon.
Global {\em a priori} $L^\infty$ bounds for $v$, as resulting from comparison principles along arbitrary smooth trajectories,
indeed mark a substantial difference between (\ref{00}) and the related chemotaxis-production systems given by
(\cite{KS}, \cite{hillen_painter})
\be{02}
	\lbal
	u_t=\na\cdot (D(u)\na u) - \na\cdot (u\na v), \\[1mm]
	\tau v_t = \Del v - f(v) + u,
	\ear
\ee
for which it is known, namely, that if $D\equiv 1$, $f\equiv id$ and $\tau\in\{0,1\}$, then already in bounded planar domains 
some solutions $(u,v)$ blow up at some finite time $T$
near which, inter alia, the spatial $L^\infty$ norm of $v$ cannot remain bounded (\cite{herrero_velazquez}, 
\cite{senba_suzuki2001}).\abs
More generally, the strongly destabilizing effect induced by the attraction-production interplay in (\ref{02}) has been reflected
in a number of results on the occurrence of blow-up in various particular scenarios. 
For the parabolic-elliptic version of (\ref{02}) with $\tau=0$ in balls $\Om\subset\R^n$ with $n\ge 2$, with radially symmetric
nonnegative initial data $u_0$, and with $f\equiv\frac{1}{|\Om|} \io u_0$, this has even led to a fairly comprehensive 
classification of nonlinearities $D$ with respect to singularity formation.
If $D$ is as in (\ref{proto}), for instance, then the exponent 
\be{03}
	\al_c^{(prod)}:=\frac{2}{n}-1 
\ee
is critical in the sense that if $\al<\frac{2}{n}-1$,
then solutions to an associated no-flux boundary value problem 
unconditionally exist globally and remain bounded, while if $\al>\frac{2}{n}-1$, then some radial initial data enforce
finite-time blow-up (\cite{cieslak_win}; cf.~also \cite{taowin_subcrit}, \cite{cieslak_stinner2012} and \cite{cieslak_stinner2015}
for some related results on a fully parabolic relative).\abs
{\bf Identifying a dimension-independent critical blow-up exponent for (\ref{00}). Main results.} \quad
Blow-up detections in chemotaxis-consumption systems so far seem limited to few results, a first class of which being concerned with
variants of (\ref{00}) that involve chemorepulsion instead of attraction (\cite{yulan_win}, \cite{ahn_win}; cf.~also
\cite{yuxiang} and \cite{dong}).
Only recently, one further finding in this direction has addressed the attraction-consumption system (\ref{00}) 
in its radial and parabolic-elliptic version, 
supplemented by boundary conditions of no-flux type for $u$ and of Dirichlet type for $v$,
thus focusing on the initial-boundary value problem
\be{0}
	\left\{ \begin{array}{ll}
	u_t = \na \cdot \big(D(u)\na u\big) - \na \cdot (u\na v),
	\qquad & x\in \Omega, \ t>0, \\[1mm]
	0 = \Delta v - uv,
	\qquad & x\in \Omega, \ t>0, \\[1mm]
	(D(u)\na u - u\na v)\cdot\nu=0, \quad v=M,
	\qquad & x\in \pO, \ t>0, \\[1mm]
	u(x,0)=u_0(x), 
	\qquad & x\in\Omega,
	\end{array} \right.
\ee
in a ball $\Om=B_R(0)\subset\R^n$ with $n\ge 1$ and $R>0$, with $M>0$, and with a given radially symmetric $u_0:\Om\to [0,\infty)$.
Indeed, a contradiction-based argument in \cite{win_bubbling} has derived a result on global nonexistence in (\ref{0}) for some radial
nonnegative $u_0\in C_0^\infty(\Om)$ 
under the assumption that $D\in C^3([0,\infty))$ is positive on $[0,\infty)$ and such that
\be{04}
	D(\xi) \le \KD (\xi+1)^{-\al}
	\qquad \mbox{for all } \xi\ge 0
\ee
with some $\KD>0$ and some
\be{05}
	\al>1.
\ee
The purpose of the present note now is to make sure that this hypothesis cannot be substantially relaxed,
and that within the framework of radially symmetric solutions, the dimension-independent number
\be{06}
	\al_c^{(cons)} := 1
\ee
plays the role of an explosion critical parameter in (\ref{0}).
The following main result in this respect embeds this into a statement on existence, uniqueness and boundedness
of global solutions for arbitrary suitably regular and radial nonnegative initial data, requiring $D$ to satisfy
a hypothesis which is essentially complementary to that expressed in (\ref{04})-(\ref{05}):
\begin{theo}\label{theo55}
  Let $n\ge 1$, $R>0$ and $\Om=B_R\subset\R^n$, let $M>0$, and assume that $D$ is such that
  \be{Dreg}
	D\in C^3([0,\infty)),
  \ee
  and that
  \be{D99}
	D(\xi)\ge \kD (\xi+1)^{-\al}
	\qquad \mbox{for all } \xi\ge 0
  \ee
  with some $\kD>0$ and some
  \bas
	\al<1.
  \eas
  Then for any $q>\max\{2,n\}$ and each nonnegative and radially symmetric $u_0\in W^{1,q}(\Om)$,
  the problem (\ref{0}) possesses a uniquely determined global classical solution $(u,v)$ with
  \bas
	\lbal
	u\in \bigcup_{q>\max\{2,n\}} C^0([0,\infty);W^{1,q}(\Om)) \cap C^{2,1}(\bom\times (0,\infty)) \qquad \mbox{and} \\[1mm]
	v \in C^{2,0}(\bom\times (0,\infty))
	\ear
  \eas
  for which $u\ge 0$ and $v\ge 0$ in $\Om\times (0,\infty)$.
  Moreover, this solution is bounded in the sense that there exists $C>0$ such that
  \be{55.1}
	\|u(\cdot,t)\|_{L^\infty(\Om)} \le C
	\qquad \mbox{for all } t>0.
  \ee
\end{theo}
\mysection{Local existence and uniqueness}
A basic but essential first step toward Theorem \ref{theo55} consists in developing an appropriate
theory of local solutions. In view of our overall purposes, a particular focus will here not only be on the derivation
of the handy extensibility criterion (\ref{ext}) below, but beyond this also on a clean argument asserting uniqueness;
in the considered parabolic-elliptic scenario, both these aspects seem seem not covered by classical literature,
neither from precedent chemotaxis analysis nor from general parabolic theory.\abs
In two essential places, our reasoning concerned with such a local theory 
will refer to the following preliminary
observation on growth of differences between solutions to the respective first sub-problems of (\ref{0}). 
\begin{lem}\label{lem99}
  Let $n\ge 1$ and $q>\max\{n,2\}$, assume (\ref{Dreg}), and let $L>0$.
  Then there exists $\Lambda(L)>0$ with the property that whenever $T>0$ and
  \bas
	\left\{ \begin{array}{l}
	\wt{u} \in C^{2,1}(\bom\times (0,T)), \\[1mm]
	\wt{v} \in C^{2,0}(\bom\times (0,T))
	\end{array} \right.
	\qquad \mbox{and} \qquad
	\left\{ \begin{array}{l}
	u \in C^0([0,T);W^{1,q}(\Om)) \cap C^{2,1}(\bom\times (0,T)), \\[1mm]
	v \in C^{2,0}(\bom\times (0,T))
	\end{array} \right.
  \eas
  are such that
  \be{99.1}
	0 \le \wt{u} \le L
	\quad \mbox{and} \quad  0\le u \le L
	\qquad \mbox{in } \Om\times (0,T)
  \ee
  as well as
  \be{99.2}
	\|\na u(\cdot,t)\|_{L^q(\Om)} \le L
	\quad \mbox{and} \quad
	\|\na v(\cdot,t)\|_{L^q(\Om)} \le L
	\qquad \mbox{for all } t\in (0,T),
  \ee
  and that
  \be{99.3}
	\left\{ \begin{array}{ll}
	\wt{u}_t = \na \cdot \big( D(\wt{u})\na \wt{u}\big) - \na \cdot (\wt{u}\na \wt{v}),
	\qquad & x\in\Om, \ t\in (0,T), \\[1mm]
	\big(D(\wt{u})\na u-\wt{u}\na \wt{v}\big)\cdot\nu=0,
	\qquad & x\in\pO, \ t\in (0,T),
	\end{array} \right.
  \ee
  as well as
  \be{99.4}
	\left\{ \begin{array}{ll}
	u_t= \na \cdot \big( D(u)\na u\big) - \na \cdot (u\na v),
	\qquad & x\in\Om, \ t\in (0,T), \\[1mm]
	\big(D(u) \na u - u\na v\big)=0,
	\qquad & x\in\pO, \ t\in (0,T),
	\end{array} \right.
  \ee
  we have
  \be{99.5}
	\frac{d}{dt} \io (\wt{u}-u)^2
	\le \Lambda(L) \io (\wt{u}-u)^2
	+ \Lambda(L) \io \big| \na (\wt{v}-v)\big|^2
	\qquad \mbox{for all } t\in (0,T).
  \ee
\end{lem}
\proof
  According to (\ref{Dreg}), given $L>0$ we can pick $c_1=c_1(L)>0$ and $c_2=c_2(L)>0$ such that
  \be{99.6}
	D(\xi)\ge c_1
	\quad \mbox{and} \quad
	|D'(\xi)| \le c_2
	\qquad \mbox{for all } \xi\in [0,L],
  \ee
  and using the inequality $\frac{2q}{q-2} < \frac{2n}{(n-2)_+}$, as asserted by our assumptions on $q$, we can draw on the
  compactness of the embedding $W^{1,2}(\Om) \hra L^\frac{2q}{q-2}(\Om)$ to find $c_3=c_3(L)>0$ fulfilling
  \be{99.7}
	\frac{(c_2^2+1) L^2}{c_1} \|\vp\|_{L^\frac{2q}{q-2}(\Om)}^2
	\le \frac{c_1}{4} \|\na\vp\|_{L^2(\Om)}^2 
	+ c_3 \|\vp\|_{L^2(\Om)}^2
	\qquad \mbox{for all } \vp\in W^{1,2}(\Om).
  \ee
  We now suppose that $T>0$, and that $(\wt{u},\wt{v})$ and $(u,v)$ satisfy (\ref{99.1}-(\ref{99.4}), and writing $w:=\wt{u}-u$
  we then infer from (\ref{99.3}) and (\ref{99.4}) that
  in $\Om\times (0,T)$,
  \bas
	w_t= \na\cdot \big(D(\wt{u})\na w\big)
	+ \na \cdot \big( [D(\wt{u})-D(u)]\na u\big)
	- \na \cdot (w\na v) - \na \cdot \big( \wt{u}\na (\wt{v}-v)\big)
  \eas
  with
  $\big(D(\wt{u})\na w + [D(\wt{u})-D(u)]\na u - w\na v -\wt{u}\na (\wt{v}-v)\big)\cdot\nu=0$ on $\pO\times (0,T)$.
  As $D(\wt{u})\ge c_1$ and $|D(\wt{u})-D(u)| \le c_2 |w|$ by (\ref{99.1}) and (\ref{99.6}), by employing Young's inequality
  we therefore obtain that 
  \bea{99.8}
	\frac{1}{2} \frac{d}{dt} \io w^2 
	+ c_1 \io |\na w|^2
	&\le& \frac{1}{2} \frac{d}{dt} \io w^2
	+ \io D(\wt{u})|\na w|^2 \nn\\
	&=& - \io \big[ D(\wt{u})-D(u)\big] \na u\cdot\na w
	+ \io w\na v\cdot\na w 
	+ \io \wt{u} \na(\wt{v}-v)\cdot\na w \nn\\
	&\le& \bigg\{ \frac{c_1}{4} \io |\na w|^2 
	+ \frac{c_2^2}{c_1} \io w^2 |\na u|^2 \bigg\}
	+ \bigg\{ \frac{c_1}{4} \io |\na w|^2
	+ \frac{1}{c_1} \io w^2 |\na v|^2 \bigg\} \nn\\
	& & + \bigg\{ \frac{c_1}{4} \io |\na w|^2
	+ \frac{L^2}{c_1} \io \big|\na (\wt{v}-v)\big|^2 \bigg\}
  \eea
  for all $t\in (0,T)$.
  Here we combine the H\"older inequality with (\ref{99.2}) and (\ref{99.7}) to estimate
  \bas
	\frac{c_2^2}{c_1} \io w^2 |\na u|^2
	+ \frac{1}{c_1} \io w^2 |\na v|^2
	&\le& \frac{c_2^2}{c_1} \|w\|_{L^\frac{2q}{q-2}(\Om)}^2 \|\na u\|_{L^q(\Om)}^2 
	+ \frac{1}{c_1} \|w\|_{L^\frac{2q}{q-2}(\Om)}^2 \|\na v\|_{L^q(\Om)}^2 \\
	&\le& \Big\{ \frac{c_2^2 L^2}{c_1} + \frac{L^2}{c_1} \Big\} \|w\|_{L^\frac{2q}{q-2}(\Om)}^2 \\
	&\le& \frac{c_1}{4} \io |\na w|^2
	+ c_3 \io w^2
  \eas
  for all $t\in (0,T)$,
  and to thereby infer (\ref{99.5}) from (\ref{99.8}) upon letting 
  $\Lambda(L):=\max \big\{ c_3(L) \, , \, \frac{L^2}{c_1(L)}\big\}$.
\qed
By an argument differing from that in classical precedents (cf., e.g., \cite{cieslak_JMAA}), 
we can now make sure that (\ref{0}) indeed allows for the following statement on local existence, uniqueness and
extensibility.
\begin{lem}\label{lem_loc}
  Assume that $n\ge 1$, that $D$ satisfies (\ref{Dreg}) and that $M>0$, and suppose that $q>\max\{2,n\}$ and
  $u_0\in W^{1,q}(\Om)$ is radially symmetric and nonnegative.
  Then there exist $\tm\in (0,\infty]$ and a pair $(u,v)$ of radial functions, for each $T\in (0,\tm]$ uniquely determined by
  the inclusions
  \be{lreg}
	\lbal
	u\in C^0([0,T);W^{1,q}(\Om)) \cap C^{2,1}(\bom\times (0,T)) \qquad \mbox{and} \\[1mm]
	v \in C^{2,0}(\bom\times (0,T)),
	\ear
  \ee
  such that $u\ge 0$ and $v\ge 0$ in $\Om\times (0,\tm)$, that $(u,v)$ solves (\ref{0}) in the classical sense in $\Om\times (0,\tm)$,
  and that
  \be{ext}
	\mbox{if $\tm<\infty$, \quad then} \quad
	\limsup_{t\nearrow\tm} \|u(\cdot,t)\|_{L^\infty(\Om)} =\infty.
  \ee
  Moreover,
  \be{mass}
	\io u(\cdot,t) dx = \io u_0 dx
	\qquad \mbox{for all } t\in (0,\tm).
  \ee
\end{lem}
\proof
  \un{Step 1.} \quad
  We first claim that for each $L>0$ there exists $T(L)>0$ with the property that whenever
  $\hnz\in W^{1,q}(\Om)$ is nonnegative with $\|\hnz\|_{L^\infty(\Om)} + \|\na\hnz\|_{L^q(\Om)} \le L$,
  the problem
  \be{l2}
	\left\{ \begin{array}{ll}
	u_t = \na \cdot \big(D(u)\na u\big) - \na \cdot (u\na v),
	\qquad & x\in \Omega, \ t\in (0,T(L)), \\[1mm]
	0 = \Del v - uv,
	\qquad & x\in \Omega, \ t\in (0,T(L)), \\[1mm]
	\big(D(u)\na u - u\na v\big)\cdot\nu=0, \quad v=M,
	\qquad & x\in \pO, \ t\in (0,T(L)), \\[1mm]
	u(x,0)=\hnz(x),
	\qquad & x\in\Omega,
	\end{array} \right.
  \ee
  admits at least one solution $0 \le u\in C^0([0,T(L)];W^{1,q}(\Om)) \cap C^{2,1}(\bom\times (0,T(L))$.\abs
  To verify this, given $L>0$ we use (\ref{Dreg}) to fix $c_i(L)>0$, $i\in\{1,2,3\}$, such that
  \be{l3}
	c_1(L) \le D(\xi) \le c_2(L)
	\quad \mbox{and} \quad
	|D'(\xi)| \le c_3(L)
	\qquad \mbox{for all } \xi\in [0,L+4],
  \ee
  and combine a Sobolev embedding theorem with elliptic regularity theory to find $c_4>0$ such that
  \be{l4}
	\|\na\psi\|_{L^\infty(\Om)}
	\le c_4 \|\Del\psi\|_{L^\infty(\Om)}
	\qquad \mbox{for all $\psi\in C^2(\bom)$ fulfilling } \psi|_{\pO}=0.
  \ee
  We then abbreviate $c_5(L):=c_4 (L+4) M$ and $c_6(L):=\frac{(L+4) c_5(L)}{c_1(L)}$ and fix some $\zeta=\zeta^{(L)}\in C^2(\bom)$
  such that $1\le \zeta\le 2$ in $\Om$ and $\frac{\pa\zeta}{\pa\nu} \ge c_6(L)$ on $\pO$.
  Letting $c_7(L):=\|\na\zeta\|_{L^\infty(\Om)}$ and $c_8(L):=\|\Del\zeta\|_{L^\infty(\Om)}$, we can thereafter choose $T=T(L)>0$
  suitably small such that
  \be{l5}
	\left\{ \begin{array}{l}
	y'(t)=\big\{ c_2(L) c_8(L) + c_5(L) c_7(L)\big\}\cdot y(t) + c_3(L) c_7^2(L) y^2(t), 	
	\qquad t\in (0,T), \\[1mm]
	y(0)=1,
	\end{array} \right.
  \ee
  possesses a solution $y\in C^1([0,T])$ fulfilling $1\le y\le 2$.
  Keeping this choice fixed, we may rely on a standard result on H\"older regularity in scalar parabolic equations
  (\cite{porzio_vespri}) to pick $\vartheta_1=\vartheta_1(L)\in (0,1)$ and $c_8(L)>0$ with the property that whenever
  $u\in C^0([0,T];W^{1,q}(\Om))\cap C^{2,1}(\bom\times (0,T))$ and $h\in C^{1,0}(\bom\times (0,T);\R^N)$ 
  are such that
  \be{l6}
	\left\{ \begin{array}{l}
	u_t = \na \cdot \big(D(u)\na u\big) - \na \cdot (u h(x,t)),
	\hs{30mm}
	x\in \Omega, \ t>0, \\[1mm]
	\big(D(u)\na u - uh(x,t)\big)\cdot\nu=0,
	\hs{42mm}
	x\in\pO, \ t\in (0,T), \\[1mm]
	\|u(\cdot,0)\|_{L^\infty(\Om)} + 
	\|\na u(\cdot,0)\|_{L^q(\Om)}
	\le L,
	\quad 
	0 \le u \le L+4
	\quad \mbox{and} \quad \\[1mm]
 	\|h(\cdot,t)\|_{L^q(\Om)} \le c_5(L) \mbox{ for all } t\in (0,T),
	\end{array} \right.
  \ee
  we have
  \be{l7}
	\|u\|_{C^{\vartheta_1,\frac{\vartheta_1}{2}}(\bom\times [0,T])} \le c_8(L).
  \ee
  We now fix any $\vartheta\in (0,\vartheta_1)$, and on the closed bounded convex subset
  \bas
	S:=\Big\{ u\in X \ \Big| \ 0 \le u \le L+4 \Big\}
  \eas
  of the Banach space $X:=C^{\vartheta,\frac{\vartheta}{2}}(\bom\times [0,T])$, we introduce a mapping $F$ by letting
  $F\wh{u}:=u$ for $\wh{u}\in S$, where according to classical parabolic theory (\cite{amann}),
  $0\le u\in C^0([0,T_1);W^{1,q}(\Om)) \cap C^{2,1}(\bom\times (0,T_1))$ is the unique solution of 
  \be{l8}
	\left\{ \begin{array}{ll}
	u_t = \na \cdot \big(D(u)\na u\big) - \na \cdot (u\na v),
	\qquad & x\in \Omega, \ t\in (0,T_1), \\[1mm]
	\big(D(u)\na u - u\na v\big)\cdot\nu=0,
	\qquad & x\in \pO, \ t\in (0,T_1), \\[1mm]
	u(x,0)=\wh{u}_0(x), 
	\qquad & x\in\Omega,
	\end{array} \right.
  \ee
  maximally extended up to some $T_1\in (0,T]$ which is such that
  \be{l9}
	\mbox{if $T_1<T$, \quad then \quad}
	\limsup_{t\nearrow T_1} \|u(\cdot,t)\|_{W^{1,q}(\Om)}=\infty,
  \ee
  where $v$ denotes the solution of the problem
  \be{l10}
	\left\{ \begin{array}{ll}
	0=\Del v - \wh{u} v,
	\qquad & x\in \Omega, \ t\in (0,T_1), \\[1mm]
	v(x,t)=M,
	\qquad & x\in \pO, \ t\in (0,T_1),
	\end{array} \right.
  \ee
  which by definition of $X$ and standard elliptic theory (\cite{GT}) indeed is solvable by some function 
  $cvin \bigcup_{\vartheta_2\in (0,1)} C^{2+\vartheta_2,\vartheta_2}(\bom\times [0,T_1))$ fulfilling 
  $0\le v \le M$ thanks to the maximum principle.
  In particular, in view of (\ref{l4}) we obtain that
  \be{l11}
	\|\na v(\cdot,t)\|_{L^\infty(\Om)} \le c_4 M(L+4) = c_5(L)
	\qquad \mbox{for all } t\in (0,T_1),
  \ee
  so that 
  \bas
	\on(x,t):=L+y(t) \zeta(x),
	\qquad x\in\bom, \ t\ge 0,
  \eas
  satisfies
  \bas
	\on_t
	- \na\cdot \big( D(\on) \na\on\big) + \na \cdot (\on\na v)
	&=& y'(t) \zeta - y(t) D(\on)\Del\zeta - y^2(t) D'(\on) |\na\zeta|^2 - y(t) \na v\cdot\na\zeta + \on \wh{n} v \\
	&\ge& y'(t) - c_2(L) c_8(L) y(t) - c_3(L) c_7^2(L) y^2(t) - c_5(L) c_7(L) y(t) \\
	&\ge& 0	\qquad \mbox{in } \Om\times (0,T_1)
  \eas
  thanks to (\ref{l5}) and the fact that 
  \be{l12}
	0\le \on \le L+4
	\qquad \mbox{in } \Om\times (0,\infty),
  \ee
  as ensured by the inequalities $0\le y\le 2$ and $0 \le \zeta\le 2$.
  Since (\ref{l11}) and (\ref{l12}) along with our restriction that $\frac{\pa\zeta}{\pa\nu}|_{\pO}\ge c_6(L)$ furthermore 
  warrant that
  \bas
	\frac{\pa\on}{\pa\nu} - \frac{\on}{D(\on)} \frac{\pa v}{\pa\nu}
	&=& y(t) \frac{\pa\zeta}{\pa\nu} - \frac{\on}{D(\on)} \frac{\pa v}{\pa \nu}  \\
	&\ge& c_6(L) y(t) - \frac{L+4}{c_1(L)} \cdot c_5(L) 
	\ge 0
	\qquad \mbox{on } \pO\times (0,T_1),
  \eas
  we can employ a comparison argument to conclude that since $\|\hnz\|_{L^\infty(\Om)} \le L$,
  we have
  \be{l13}
	u \le \on
	\qquad \mbox{in } \Om\times (0,T_1).
  \ee
  This especially means that the parabolic problem (\ref{l8}) actually is non-degenerate, so that standard parabolic theory
  applies so as to guarantee that $u\in C^{2,1}(\bom\times [\frac{T_1}{2},T_1])$, and that thus, by (\ref{l9}),
  in fact we must have $T_1=T$.\abs
  As (\ref{l13}) and (\ref{l12}) thus imply that $0\le u\le L+4$ in $\Om\times (0,T)$, in view of (\ref{l8}) and (\ref{l6})
  we may now draw on (\ref{l7}) to conclude that $F$ indeed maps $S$ into itself, and that moreover $F(S)$ is bounded
  in $C^{\vartheta_1,\frac{\vartheta_1}{2}}(\bom\times [0,T])$, thanks to the inequality $\vartheta_1>\vartheta$ meaning that
  $\ov{F(S)}$ is a compact subset of $X$.\abs
  To see that $F$ is continuous, we assume that $(\wh{u}_j)_{j\in\N} \subset S$ and $\wh{u}\in X$ be such that
  \be{l110}
	\wh{u}_j \to \wh{u}
	\quad \mbox{in } X
	\qquad \mbox{as } j\to\infty,
  \ee
  and letting $u_j:=F(\wh{u}_j)$ for $j\in\N$ we first observe that,
  again by elliptic regularity theory, there exists $\vartheta_2\in (0,1)$ such that for all $j\in\N$ one can find
  $v_j\in C^{2+\vartheta_2,\vartheta_2}(\bom\times [0,T])$ such that $0=\Del v_j - \wh{u}_j v_j$ in $\Om\times (0,T)$ and
  $v_j|_{\pO}=M$, and that $v_j\to v$ in $C^{2+\vartheta_2,\vartheta_2}(\bom\times [0,T])$ as $j\to\infty$, with some
  classical solution $v$ of $0=\Del v - \wh{u} v$ in $\Om\times (0,T)$ fulfilling $v|_{\pO}=M$.
  On the other hand, using these regularity features of $v$ we may once again rely on standard parabolic theory (\cite{amann})
  to infer the existence of a classcial solution 
  $u\in C^0([0,T];W^{1,q}(\Om)) \cap C^{2,1}(\bom\times (0,T))$ of the corresponding version of (\ref{l8}) with $T_1=T$,
  by definition of $F$ meaning that we must have $u=F(\wh{u})$.
  In particular, besides $\sup_{t\in (0,T)} \|\na v_j(\cdot,t)\|_{L^q(\Om)}$, also
  both $\sup_{t\in (0,T)} \|u(\cdot,t)\|_{L^\infty(\Om)}$ and $\sup_{t\in (0,T)} \|\na u(\cdot,t)\|_{L^q(\Om)}$
  are finite, so that since due to the inclusion $F(S)\subset S$ we already know that $0\le u_j \le L+4$ for all $j\in\N$,
  we may employ Lemma \ref{lem99} to find $c_9>0$ such that
  \be{l112}
	\frac{d}{dt} \io (u_j-u)^2 \le c_9 \io (u_j-u)^2 + c_9 \io \big|\na (v_j-v)\big|^2
	\qquad \mbox{for all } t\in (0,T).
  \ee
  But since $v_j-v=0$ on $\pO\times (0,T)$, a Poincar\'e inequality in $W_0^{1,2}(\Om)$ becomes applicable
  so as to assert that when tesing the identity
  \bas
	-\Del (v_j-v)=-(\wh{u}_j-\wh{u}) v - \wh{u}_j (v_j-v)
  \eas
  against $v_j-v$ and using that $0\le v\le M$ and $\wh{u}_j\ge 0$, with some $c_{10}>0$ we see that
  \bea{l113}
	\io \big|\na (v_j-v)\big|^2
	&=& - \io (\wh{u}_j-\wh{u}) v(v_j-v)
	- \io \wh{u}_j (v_j-v)^2 \nn\\
	&\le& M \|\wh{u}_j-u\|_{L^2(\Om)} \|v_j-v\|_{L^2(\Om)} \nn\\
	&\le& c_{10} \|\wh{u}_j-u\|_{L^2(\Om)} \|\na(v_j-v)\|_{L^2(\Om)} \nn\\
	&\le& \frac{1}{2} \io \big|\na (v_j-v)\big|^2
	+ \frac{c_{10}^2}{2} \io (\wh{u}_j-\wh{u})^2
  \eea
  for all $t\in (0,T)$.
  From (\ref{l112}) we therefore obtain that 
  \bas	
	\frac{d}{dt} \io (u_j-u)^2 \le c_9 \io (u_j-u)^2
	+ \frac{c_9 c_{10}^2}{2} \io (\wh{u}_j-\wh{u})^2
  \eas
  for all $t\in (0,T)$ and $j\in\N$,
  from which in view of the identity $(u_j-u)|_{t=0}=0$ and (\ref{l110}) we readily infer that as $j\to\infty$ we have
  $F(\wh{u}_j)=u_j\to u=F(\wh{u})$ in 
  $L^\infty((0,T);L^2(\Om))$ and hence, by relative compactness of $(u_j)_{j\in\N} \subset F(S)$ in $X$, also
  in $X$.\abs
  \un{Step 2.} \quad
  Let us next make sure that given any nonnegative $u_0\in W^{1,q}(\Om)$
  we can find $\tm\in (0,\infty]$ and a classical
  solution $(u,v)$ of (\ref{0}) in $\Om\times (0,\tm)$ which satisfies (\ref{lreg}) with $u\ge 0$, and for which
  (\ref{ext}) holds.\abs
  Indeed, a standard extension procedure on the basis of Step 1 provides $\tm\in (0,\infty]$ and a classical solution
  $(u,v)$ on $\Om\times (0,\tm)$ fulfilling (\ref{lreg}) and $u\ge 0$ as well as the alternative that
  \be{l23}
	\mbox{if $\tm<\infty$, \quad then} \quad
	\limsup \|u(\cdot,t)\|_{W^{1,q}(\Om)} =\infty.
  \ee
  To see that actually (\ref{ext}) holds, assuming on the contrary that $\tm$ 
  and $\sup_{t\in (0,\tm)} \|u(\cdot,t)\|_{L^\infty(\Om)}$ both be finite, we could 
  again invoke arguments from elliptic and parabolic regularity theories to successively infer that
  $\na v\in L^\infty(\Om\times (\frac{1}{16}\tm,\tm))$ (\cite{GT}), that
  $u\in \bigcup_{\vartheta\in (0,1)} C^{\vartheta,\frac{\vartheta}{2}}(\bom\times [\frac{1}{8}\tm,\tm])$ (\cite{porzio_vespri}),
  that $v\in \bigcup_{\vartheta\in (0,1)} C^{2+\vartheta,\vartheta}(\bom\times [\frac{1}{4} \tm,\tm])$ (\cite{GT}), and that
  $u\in \bigcup_{\vartheta\in (0,1)} C^{1+\vartheta,\frac{1+\vartheta}{2}}(\bom\times [\frac{1}{2} \tm,\tm])$ (\cite{lieberman}).
  The latter would contradict (\ref{l23}), however.\abs
  \un{Step 3.} \quad
  It remains to verify the claimed uniqueness feature, as then clearly also the statement on radial symmetry follows.\abs
  To this end, supposing that for some $T>0$ we are given two solutions $(u_i,v_i)$ of (\ref{0}) in $\Om\times (0,T)$
  fulfilling (\ref{lreg}),
  $i\in \{1,2\}$, we may apply Lemma \ref{lem99} for a second time to see that for each $T_0\in (0,T)$ we can find
  $c_{11}=c_{11}(T_0)>0$ such that $w:=u_1-u_2$ satisfies
  \be{l33}
	\frac{d}{dt} \io w^2 \le c_{11} \io w^2 + c_{11} \io \big|\na (v_1-v_2)\big|^2
	\qquad \mbox{for all } t\in (0,T_0).
  \ee 
  But since an almost verbatim copy of the reasoning in (\ref{l113}) shows that with some $c_{12}>0$ we have
  \bas
	\io \big| \na (v_1-v_2)\big|^2 \le c_{12} \io w^2
	\qquad \mbox{for all } t\in (0,T),
  \eas
  from (\ref{l33}) we infer that $\frac{d}{dt} \io w^2 \le c_{11}(1+c_{12}) \io w^2$ for all $t\in (0,T_0)$,
  upon an integration using that $w|_{t=0}=0$ showing that, in fact, $w\equiv 0$ in $\Om\times (0,T_0)$ for any such $T_0$.
\qed
\mysection{Proof of Theorem \ref{theo55}}
As is quite common in the analysis of quasilinear chemotaxis systems, our reasoning toward global extensibility of
the local solution found above will be arranged around the derivation of time-independent $L^p$ bounds for its first
component, with $p>1$ being arbitrarily large but finite.
Our first step in this direction already uses radial symmetry in order to infer a preliminary upper estimate
for a boundary integral that appears as the only ill-signed summand in an identity obtained by means of a 
standard testing procedure.
At this stage, the hypotheses of Theorem \ref{theo55} are not used in their full strength yet, however;
indeed, the following does not impose any restrictions on the number $\al$ in (\ref{D99}).
\begin{lem}\label{lem1}
  Let $n\ge 1$ and $M>0$, assume (\ref{Dreg}) and (\ref{D99}) with some $\al\in\R$ and $\kD>0$, 
  let $q>\max\{2,n\}$ and $u_0\in W^{1,q}(\Om)$ be radially symmetric and nonnegative,
  and let $\tm$ and $(u,v)$ be as in Lemma \ref{lem_loc}. 
  Then for each $p>1$ there exists $C(p)>0$ with the property that
  \be{1.1}
	\frac{d}{dt} \io (u+1)^p dx 
	+ \frac{1}{C(p)} \io \big|\na (u+1)^\frac{p-\al}{2}\big|^2 dx 
	\le C(p) \int_{\pO} (u+1)^p dx
	\qquad \mbox{for all } t\in (0,\tm).
  \ee
\end{lem}
\proof
  An integration by parts using (\ref{0}) shows that if given $p>1$ we let
  \be{1.2}
	\psi(\xi):=\int_0^\xi \sig (\sig+1)^{p-2} d\sig,
	\qquad \xi\ge 0,
  \ee
  then
  \bea{1.3}
	\frac{d}{dt} \io (u+1)^p dx 
	&=& p \io (u+1)^{p-1} \na \cdot \big\{ D(u)\na u - u\na v\big\} dx \nn\\
	&=& - p(p-1) \io (u+1)^{p-2} D(u) |\na u|^2 dx
	+ p(p-1) \io u(u+1)^{p-2} \na u\cdot\na v dx \nn\\
	&=& - p(p-1) \io (u+1)^{p-2} D(u) |\na u|^2 dx
	+ p(p-1) \io \na \psi(u)\cdot\na v dx \nn\\
	&=& - p(p-1) \io (u+1)^{p-2} D(u) |\na u|^2 dx
	- p(p-1) \io \psi(u) \Del v dx \nn\\
	& & + p(p-1) \int_{\pO} \psi(u) \frac{\pa v}{\pa\nu} dx \nn\\
	&=& - p(p-1) \io (u+1)^{p-2} D(u) |\na u|^2 dx
	- p(p-1) \io \psi(u) \cdot uv dx \nn\\
	& & + p(p-1) \int_{\pO} \psi(u) \frac{\pa v}{\pa\nu} dx 
  \eea
  for all $t\in (0,\tm)$.
  In order to estimate the rightmost summand herein, 
  we first observe that an integration of the identity 
  $\pa_r (r^{n-1} v_r) = r^{n-1} uv$,
  valid for $r\in (0,R)$ and $t\in (0,\tm)$ according to the second equation in (\ref{0}), shows that
  \bas
	v_r(r,t)=r^{1-n} \int_0^r \rho^{n-1} u(\rho,t)v(\rho,t) d\rho
	\qquad \mbox{for all $r\in (0,R]$ and } t\in (0,\tm).
  \eas
  Using that $0\le v\le M$ by comparison, in line with (\ref{mass}) we thus infer on writing 
  $c_1:=\frac{R^{1-n}}{n|B_1(0)|} \io u_0 dx$ that
  \bas
  	\frac{\pa v}{\pa\nu}\le c_1
	\qquad \mbox{on } \pO\times (0,\tm).
  \eas
  As furthermore, by (\ref{1.2}),
  \bas
	0 \le \psi(\xi)
	&\le& \int_0^\xi (\sig+1)^{p-1} d\sig \\
	&\le& \frac{1}{p}(\xi+1)^p
	\qquad \mbox{for all } \xi\ge 0,
  \eas
  we can estimate
  \bas
	p(p-1) \int_{\pO} \psi(u) \frac{\pa v}{\pa\nu} dx
	&\le& (p-1) c_1 \int_{\pO} (u+1)^p dx
	\qquad \mbox{for all } t\in (0,\tm),
  \eas
  so that (\ref{1.1}) becomes a consequence of (\ref{1.3}) and the fact that
  \bas
	p(p-1) \io (u+1)^{p-2} D(u) |\na u|^2 dx
	&\ge&
	p(p-1) \kD \io (u+1)^{p-\al-2} |\na u|^2 dx \nn\\
	&=& \frac{4p(p-1)\kD}{(p-\al)^2} \io \big| \na (u+1)^\frac{p-\al}{2}\big|^2
	\qquad \mbox{for all } t\in (0,\tm)
  \eas
  thanks to (\ref{D99}),
  because clearly
  \bas
	- p(p-1) \io \psi(u) \cdot uv dx 
	\le 0
	\qquad \mbox{for all } t\in (0,\tm)
  \eas
  by nonnegativity of $\psi$.
\qed
Suitably controlling the boundary integral on the right-hand side in (\ref{1.1}) by means of the dissipated quantity therein now 
requires the full set of assumptions made in Theorem \ref{theo55}:
In the range of $\al$ addressed there, namely, we may once more draw on radial symmetry to estimate the boundary values
of $u$ by using a Gagliardo-Nirenberg interpolation in the one-dimensional interval $(\frac{R}{2},R)$,
and thereby achieve the following.
\begin{lem}\label{lem2}
  If $n\ge 1$ and $M>0$, if (\ref{Dreg}) and (\ref{D99}) hold with some 
  \be{2.01}
	\al<1
  \ee
  and some $\kD>0$,
  and if $q>\max\{2,n\}$ and $u_0\in W^{1,q}(\Om)$ is radially symmetric and nonnegative,
  then for every $p>1$ any any $\eta>0$
  there exists $C(p,\eta)>0$ such that with $\tm$ and $(u,v)$ taken from Lemma \ref{lem_loc}, we have
  \be{2.1}
	\int_{\pO} (u+1)^p dx
	\le \eta \io \big|\na (u+1)^\frac{p-\al}{2}\big|^2 dx 
	+ C(p,\eta)
	\qquad \mbox{for all } t\in (0,\tm).
  \ee
\end{lem}
\proof
  Thanks to radial symmetry, we may draw on the one-dimensional Gagliardo-Nirenberg inequality to see that due to Young's 
  inequality, with some $c_1=c_1(p)>0$ and $c_2=c_2(p)>0$ we have
  \bea{55.3}
	\int_{\pO} (u+1)^p dx
	&\le& |\pO| \cdot \|u+1\|_{L^\infty(\pO)} \nn\\
	&\le& |\pO| \cdot \|u+1\|_{L^\infty(\Om\sm B_\frac{R}{2})} \nn\\
	&=& \big\| (u+1)^\frac{p-\al}{2} \big\|_{L^\infty(\Om\sm B_\frac{R}{2})}^\frac{2p}{p-\al} \nn\\
	&\le& c_1 \big\| \pa_r (u+1)^\frac{p-\al}{2} \big\|_{L^2((\frac{R}{2},R))}^\frac{2p}{p+1-\al}
		\big\| (u+1)^\frac{p-\al}{2} \big\|_{L^\frac{2}{p-\al}((\frac{R}{2},R))}^\frac{2p}{(p-\al)(p+1-\al)} \nn\\
	&+& c_1 \big\| (u+1)^\frac{p-\al}{2} \big\|_{L^\frac{2}{p-\al}((\frac{R}{2},R))}^\frac{2p}{p-\al}  \nn\\
	&\le& c_2 \big\|\pa_r (u+1)^\frac{p-\al}{2} \big\|_{L^2((\frac{R}{2},R))}^\frac{2p}{p+1-\al} + c_2
	\qquad \mbox{for all } t\in (0,\tm),
  \eea
  because $p>1>\al$, and because
  \bas
	\big\| (u+1)^\frac{p-\al}{2} \big\|_{L^\frac{2}{p-\al}((\frac{R}{2},R))}^\frac{2}{p-\al}
	&=& \int_\frac{R}{2}^R (u+1) dr \\
	&\le& \Big(\frac{2}{R}\Big)^{n-1} \int_\frac{R}{2}^R r^{n-1} (u+1) dr  \\
	&\le& \frac{(\frac{2}{R})^{n-1}}{n|B_1(0)|} \cdot \io (u+1) dx \\
	&=& \frac{(\frac{2}{R})^{n-1}}{n|B_1(0)|} \cdot \io (u+1) dx
	\qquad \mbox{for all } t\in (0,\tm)
  \eas
  by (\ref{mass}).
  We now make full use of the hypothesis in (\ref{2.01}) in observing that $\frac{2p}{p+1-\al} < 2$, so that we may employ
  Young's inequality to see that given $\eta>0$ we can find $c_3=c_3(p,\eta)>0$ fulfilling
  \bas
	c_2 \big\|\pa_r (u+1)^\frac{p-\al}{2} \big\|_{L^2((\frac{R}{2},R))}^\frac{2p}{p+1-\al}
	&\le& \frac{n|B_1(0)|}{(\frac{2}{R})^{n-1}} \cdot \eta \cdot \big\|\pa_r (u+1)^\frac{p-\al}{2} \big\|_{L^2((\frac{R}{2},R))}^2 
	+ c_3 \\
	&=& \frac{n|B_1(0)|}{(\frac{2}{R})^{n-1}} \cdot \eta \cdot \int_{\frac{R}{2}}^R \big| \pa_r (u+1)^\frac{p-\al}{2} \big|^2 dr
	+ c_3 \\
	&\le& n|B_1(0)| \cdot \eta \cdot \int_{\frac{R}{2}}^R r^{n-1} \big| \pa_r (u+1)^\frac{p-\al}{2} \big|^2 dr
	+ c_3 \\
	&\le& \eta \io \big| \na (u+1)^\frac{p-\al}{2} \big|^2 dx
	+ c_3 
	\qquad \mbox{for all } t\in (0,\tm).
  \eas
  Consequently, (\ref{2.1}) results from (\ref{55.3}).
\qed
To finally turn (\ref{1.1}) into an autonomous ODI for $t\mapsto \io (u+1)^p dx$, we interpolate on the full domain $\Om$,
and may here proceed in a comparatively less subtle manner, inter alia not requiring the inequality $\al<1$ in this step.
\begin{lem}\label{lem3}
  Let $n\ge 1$ and $M>0$, suppose that (\ref{Dreg}) and (\ref{D99}) hold with some $\al\in\R$ and some $\kD>0$, and
  let $q>\max\{2,n\}$ and $u_0\in W^{1,q}(\Om)$ be radially symmetric and nonnegative.
  Then given any $p>1$ such that $p>\frac{n\al}{2}$, one can find $C(p)>0$ with the property  
  that the quantities $\tm$ and $u$ from Lemma \ref{lem_loc} satisfy
  \be{3.1}
	\bigg\{ \io (u+1)^p dx \bigg\}^\frac{2+n(p-1-\al)}{n(p-1)}
	\le C(p) \io \big| \na (u+1)^\frac{p-\al}{2} \big|^2 dx + C(p)
	\qquad \mbox{for all } t\in (0,\tm).
  \ee
\end{lem}
\proof
  This follows upon observing that since the inequality $p>\frac{n\al}{2}$ warrants that $\frac{2p}{p-\al}<\frac{2n}{(n-2)_+}$,
  we may combine the Gagliardo-Nirenberg inequality in its $n$-dimensional version with (\ref{mass}) to find $c_1=c_1(p)>0$ and 
  $c_2=c_2(p)>0$ fulfilling
  \bas
	\bigg\{ \io (u+1)^p dx \bigg\}^\frac{2+n(p-1-\al)}{n(p-1)}
	&=& \big\| (u+1)^\frac{p-\al}{2} \big\|_{L^\frac{2p}{p-\al}(\Om)}^{\frac{2p}{p-\al} \cdot \frac{2+n(p-1-\al)}{n(p-1)}} \\
	&\le& c_6 \big\| \na (u+1)^\frac{p-\al}{2} \big\|_{L^2(\Om)}^2
		\big\| (u+1)^\frac{p-\al}{2} \big\|_{L^\frac{2}{p-\al}(\Om)}^\frac{2(2p-n\al)}{n(p-\al)(p-1)} \\
	& & + c_6 \big\| (u+1)^\frac{p-\al}{2} 
		\big\|_{L^\frac{2}{p-\al}(\Om)}^{\frac{2p}{p-\al}\cdot \frac{2+n(p-1-\al)}{n(p-1)}} \\
	&\le& c_7 \io \big| \na (u+1)^\frac{p-\al}{2} \big|^2 dx + c_7
  \eas
  for all $t\in (0,\tm)$.
\qed
It remains to combine the previous three lemmata with a standard argument turning $L^p$ estimates into $L^\infty$ bounds
by a Moser iteration, leading to the announced result on global existence and boundedness of radial solutions to
(\ref{0}) throughout the considered range of subcritical $\al$:\abs
\proofc of Theorem \ref{theo55}. \quad
  Taking $\tm$ and $(u,v)$ from Lemma \ref{lem_loc}, in view of (\ref{ext}) we only need to make sure that
  \be{55.11}
	\sup_{t\in (0,\tm)} \|u(\cdot,t)\|_{L^\infty(\Om)} <\infty.
  \ee
  To this end, we first fix an arbitrary $p>1$ such that $p>\frac{n\al}{2}$, and then obtain from Lemma \ref{lem1}
  that with some $c_1=c_1(p)>0$ and $c_2=c_2(p)>0$ such that
  \be{55.61}
	\frac{d}{dt} \io (u+1)^p dx 
	+ c_1 \io \big|\na (u+1)^\frac{p-\al}{2}\big|^2 dx 
	\le c_2 \int_{\pO} (u+1)^p dx
	\qquad \mbox{for all } t\in (0,\tm),
  \ee
  whereupon an application of Lemma \ref{lem2} yields $c_3=c_3(p)>0$ fulfilling
  \be{55.62}
	c_2 \int_{\pO} (u+1)^p dx
	\le \frac{c_1}{2} \io \big|\na (u+1)^\frac{p-\al}{2}\big|^2 dx 
	+ c_3
	\qquad \mbox{for all } t\in (0,\tm).
  \ee
  Relying on the inequality $p>\frac{n\al}{2}$, from Lemma \ref{lem3} we then infer the existence of $c_4=c_4(p)>0$ such that
  \be{55.63}
	\frac{c_1}{2} \io \big|\na (u+1)^\frac{p-\al}{2}\big|^2 dx 
	\ge c_4 \cdot \bigg\{ \io (u+1)^p dx \bigg\}^\frac{2+n(p-1-\al)}{n(p-1)}
	- \frac{c_1}{2}
	\qquad \mbox{for all } t\in (0,\tm),
  \ee
  whence collecting (\ref{55.61})-(\ref{55.63}) we see that
  \bas
	\frac{d}{dt} \io (u+1)^p dx 
	+ c_4 \cdot \bigg\{ \io (u+1)^p dx \bigg\}^\frac{2+n(p-1-\al)}{n(p-1)}
	\le c_3+\frac{c_1}{2}
	\qquad \mbox{for all } t\in (0,\tm),
  \eas
  from which through an ODE comparison argument it follows that
  \be{55.5}
	\sup_{t\in (0,\tm)} \io (u+1)^p dx <\infty.
  \ee
  Applying this first to an arbitrary fixed $p>n$, from this we conclude
  by means of elliptic regularity theory and thanks to the continuity of the
  embedding $W^{2,p}(\Om)\hra W^{1,\infty}(\Om)$ that also 
  \bas
	\sup\limits_{t\in (0,\tm)} \|\na v(\cdot,t)\|_{L^\infty(\Om)} < \infty,
  \eas
  whereupon using (\ref{55.5}) for a second time we may rely on the outcome of a Moser-type iterative reasoning
  (cf.~\cite{ding_win} for a statement covering the present situation) that indeed (\ref{55.11}) holds.
\qed

\bigskip

{\bf Acknowledgement.} \quad
  The author acknowledges support of the {\em Deutsche Forschungsgemeinschaft} (Project No.~462888149).
  He moreover declares that he has no conflict of interest.\abs
{\bf Data availability statement.} \quad
Data sharing is not applicable to this article, as no datasets were
generated or analyzed during the current study.\abs

\end{document}